\documentclass[12pt]{article}
\usepackage{amsmath}
\usepackage{latexsym}
\usepackage{amssymb}

\def\<{{\langle}}
\def\>{{\rangle}}
\def\pf{{\it Proof. }}
\def\C{{\,\subseteq\,}}

\def\map{{\longrightarrow}}
\def\cprime{$'$}

\def\Alg{{\mbox{\rm Alg}}}
\def\Aut{{\mbox{\rm Aut}}}
\def\End{{\mbox{\rm End}}}

\def\lcm{{\mbox{\rm lcm}}}

\def\qed{{\mbox{\ \ \ \ \ \ \ $\blacksquare$} }}

\def\L{{\Lambda}}

\def\e{{\varepsilon}}
\def\a{{\alpha}}
\def\b{{\beta}}
\def\d>{{\rightharpoondown}}
\def\u>{{\rightharpoonup}}
\def\<{{\leftharpoondown}}

\def\l{{\lambda}}

\def\L{{\Lambda}}

\def\w{{\omega}}
\def\ba{{\mathbf a}}
\def\bb{{\mathbf b}}
\def\bx{{\mathbf x}}

\def\hbar{{\overline{h}}}

\def\K{{\mathcal K}}

\def\Q{{\bf Q}}
\def\Z{{\bf Z}}

\def\tr{{\mbox{\rm Tr}}}

\numberwithin{equation}{section}

\newtheorem{thm}{Theorem}[section]

\newtheorem{prop}[thm]{Proposition}
\newtheorem{remark}[thm]{Remark}
\newtheorem{lem}[thm]{Lemma}

\newtheorem{example}[thm]{Example}

\begin{document}

\author{Siu-Hung Ng \\ Mathematics Department, Towson University, Baltimore,
MD 21252}
\title{Non-semisimple Hopf Algebras of Dimension $p^2$}
\date{}
 \maketitle
\begin{abstract}
Let $H$ be a Hopf algebra of dimension $pq$ over an algebraically
closed field of characteristic 0, where $p \le q$ are odd primes.
Suppose that $S$ is the antipode of $H$. If $H$ is not semisimple,
then $S^{4p}=id_H$ and $\tr(S^{2p})$ is an integer divisible by
$p^2$. In particular, if  $\dim H = p^2$, we prove that $H$ is
isomorphic to a Taft algebra. We then complete the classification
for the Hopf algebras of dimension $p^2$.
\end{abstract}

\section{Introduction} Let $p$ be a prime number and $k$ an
algebraically closed field of charactericstic 0. If $H$ is a
semisimple Hopf algebra of dimension $p^2$, then $H$ is isomorphic
to a group algebra \cite{Mas96}, namely $k[\Z_{p^2}]$ or $k[\Z_p
\times \Z_p]$. For the Hopf algebras $H$ of dimension $p^2$, the
only known non-semisimple Hopf algebras of dimension $p^2$ are the
Taft algebras \cite{Taft71} (cf. \cite[5]{Mont98}). The question
whether the Taft algebras are the only non-semisimple Hopf
algebras of dimension $p^2$ is open. In fact, it is also a
question suggested by Susan Montgomery in several international
conferences. It was proved in \cite[Theorem A]{AnSc98} that if
both $H$ and $H^*$ have nontrivial group-like elements or the order of
the antipode is $2p$, then $H$ is isomorphic to a Taft algebra
provided $\dim H=p^2$. In this paper, we will give a complete
answer to the question. More explicitly, we prove that for any
non-semisimple Hopf algebra $H$ over $k$
 of dimension $p^2$, $H$ is isomorphic to a Taft algebra. Hence,
 the Hopf algebras over $k$ of dimension $p^2$ can be completely
 classified (Theorem 6.5).\\

If $p \le q$ are odd primes, whether there is a non-semisimple
Hopf algebra of dimension $pq$ other than the Taft algebras is
still in question. Nevertheless, we prove for any Hopf algebra of
this type, the order of its antipode $S$ divides $4p$. Moreover,
$\tr(S^{2p})$ is an integer divisible by $p^2$ (Theorem 6.4). The
uniqueness of
Taft algebras is a consequence of this result.\\

The article is organized as follows: In section \ref{s1}, we
recall some notation, general theorems and  some useful
statements. In section \ref{s2}, we introduce the notion of the
{\em index} of a Hopf algebra and we compute the index of the Taft
algebras. In section \ref{s3}, we consider the common eigenspaces
of $S^2$ and $r(g)$ where  $S$ and $r(g)$ are the antipode and the
right multiplication by the distinguished group-like element $g$
of the Hopf algebra $H$.  We derive some arithmetic properties of
the dimensions of these eigenspaces for the Hopf algebras of odd
index. We further exploit the arithmetic properties of these
numbers
 for Hopf algebras of odd prime index in section
\ref{s4}. Finally, we prove our main theorems in
section \ref{s5}.\\

\section{Notation and Preliminaries}\label{s1}
Throughout this paper $k$ is an algebraically closed field of
characteristic 0 and $H$ is a finite-dimensional Hopf algebra over
$k$ with antipode $S$. Its comultiplication and counit are,
respectively, denoted by $\Delta$ and  $\e$. We will use
Sweedler's notation \cite{Sw69}:
$$
\Delta(x) = \sum x_{(1)} \otimes x_{(2)}\,.
$$
A non-zero element $a \in H$ is called group-like if $\Delta(a)=a
\otimes a$.  For the details of elementary aspects for
finite-dimensional Hopf algebras,  readers are referred to
the references \cite{Sw69} and \cite{Mont93bk}.\\

The set of all group-like elements $G(H)$ of $H$ is a linearly
independent set and forms a group under the multiplication of $H$.
The divisibility of $\dim H$ by $|G(H)|$ is an immediate
consequence of the following generalization of Lagrange's theorem,
due to Nichols and Zoeller:
\begin{thm}{\rm \cite{Nich89}}\label{Nich89}
If $B$ is a Hopf subalgebra of $H$, then $H$ is a free $B$-module.
In particular, $\dim B$ divides $\dim H$.
\end{thm}

The order of the antipode is of fundamental importance to the
semisimplicity of $H$. We recall some important results on the
antipode $S$ of finite-dimensional Hopf algebras $H$.

\begin{thm}{\rm \cite{LaRa87, LaRa88}} \label{th1.1}
Let $H$ be a finite dimensional Hopf algebra over a field of
characteristic 0. Then the following are equivalent:
\begin{enumerate}
\item[\rm (i)] $H$ is semisimple.
\item[\rm (ii)] $H^*$ is semisimple.
\item[\rm (iii)] $\tr(S^2) \ne 0$\,.
\item[\rm (iv)] $S^2 =id_H$\,.\\
\end{enumerate}
\end{thm}

 Let $\l$ be a non-zero right integral of $H$ and let $\L$ be a non-zero
 left integral of $H^*$. There is an $\a \in \Alg(H,k)=G(H^*)$, independent of
the choice of $\L$,  such that $\L a = \a(a) \L$ for  $a \in H$.
Likewise, there is a group-like element $g \in H$, independent of
the choice of $\l$,  such that $\b\l  = \b(g)\l$ for  $\b\in H^*$.
We call $g$ the distinguished group-like element of $H$ and $\a$
the distinguished group-like element of $H^*$. Then we have a
formula for $S^4$ in terms of $\a$ and $g$~\cite{Radf76}:
\begin{equation}\label{eS}
S^4(h) = g(\a \rightharpoonup h \leftharpoonup \a^{-1}) g^{-1}
\quad\mbox{for }a \in H\,,
\end{equation}
where $\rightharpoonup$ and $\leftharpoonup$ denote the natural
actions Hopf algebra $H^*$ on $H$ described by
$$
\b\rightharpoonup a = \sum a_{(1)}\b(a_{(2)}) \quad \mbox{and}
\quad a \leftharpoonup \b = \sum \b(a_{(1)})a_{(2)}
$$
for $\b \in H^*$ and $a \in H$. If $\l$ and $\L$ are normalized,
there are formulae for the trace of any linear endomorphism on
$H$.
\begin{thm}{\rm \cite[Theorem 2]{Radf94}}\label{th1.2}
Suppose that $\l(\L)=1$. Then for any $f \in \End_k(H)$,
\begin{eqnarray*}
\tr(f)&=&\sum \l\left(S(\L_{(2)})f(\L_{(1)}\right))\\
        &=&\sum \l\left(S\circ f(\L_{(2)})\L_{(1)}\right)\\
        &=&\sum \l\left(f\circ S(\L_{(2)})\L_{(1)}\right)\,.
\end{eqnarray*}
\end{thm}

We shall also need the following lemma of linear algebra:
\begin{lem}{\rm \cite[Lemma 2.6] {AnSc98}}\label{1}
Let $T$ be an operator on a finite dimensional vector space $V$
over $k$. Let $p$ be an odd prime and let $\w \in k$ be a
primitive $p$th root of unity.
\begin{enumerate}
\item[\rm (i)]If $\tr(T)=0$ and $T^p=id_V$, then  $\dim V_i$ is constant
where $V_i$ is eigenspace of $T$ associated with the eigenvalue
$\w^i$. In particular, $p | \dim V$.
\item[\rm (ii)]
If $\tr(T)=0$ and $T^{2p}=id_V$, then
$$
\tr(T^p) = pd
$$
for some integer $d$.
\end{enumerate}
\end{lem}
\section{Index of a Hopf algebra}\label{s2}
The distinguished group-like element $g$ defines a coalgebra
automorphism $r(g)$ on $H$ as follows:
$$
r(g)(a) = ag \quad\mbox{for}\quad a \in H\,.
$$
Since $S^2$ is an algebra automorphism on $H$,
$$
S^2 \circ r(g) = r(g)\circ S^2\,.
$$
Moreover, both $S^2$ and $r(g)$ are of finite order. Therefore,
$S^2$ and $r(g)$ generate a finite abelian subgroup of
$\Aut_k(H)$. We will simply call the exponent of the subgroup
generated by $S^4$ and $r(g)$ the {\em index} of $H$. It is easy
to see that the index of $H$ is also the smallest positive integer
$n$ such that
$$
S^{4n} = id_H \quad \mbox{and}\quad g^n=1\,.
$$
Obviously, $o(g)\mid n$ and $o(S^4)\mid n$, where $o(g)$ and $
o(S^4)$ are the orders of $g$ and $S^4$ respectively. By equation
(\ref{eS}),
\begin{equation}\label{eqdiv}
n\mid\lcm(o(g),o(\a))\,.
\end{equation}
\begin{example}\label{index_ex}\hspace{1pt}
{\rm
\begin{enumerate}
\item[\rm (i)] If both $H$ and $H^*$ are unimodular,
then $S^4=id_H$ by (\ref{eS}). Therefore, the index of $H$ is 1.
In particular, if $H$ is semisimple, the index of $H$ is 1.
\item[\rm (ii)] Let $\xi\in k$ be a $n$th root of unity. The
Taft algebra~\cite{Taft71}  $T(\xi)$ over $k$ is generated by $x$
and $a$, as a $k$-algebra, subject to the relations
$$
a^n=1,\quad ax=\xi xa,\quad x^n=0\,.
$$
The Hopf algebra structure is given by
\begin{eqnarray*}
\Delta(a)=a \otimes a,& S(a)=a^{-1}, & \e(a)=1, \\
\Delta(x) = x \otimes a +1 \otimes x, & S(x)=-xa^{-1},& \e(x)=0\,.
\end{eqnarray*}
It is known that $\{x^ia^j\,|\,0\le i,j \le n-1\}$ forms a basis
for $T(\xi)$. In particular, $\dim T(\xi) =n^2$. The linear
functional $\l$, defined by
$$
\l(x^i a^j) =\delta_{i,n-1}\delta_{j,0}
$$
is a right integral for $T(\xi)$. One can easily see that $a$ is
the distinguished group-like element of $T(\xi)$. Moreover,
$S^4(x)=\xi^2 x$ and $S^4(a)=a$. Therefore, the order of $S^4$ is
$\frac{n}{\gcd(2,n)}$. Since the order of $a$ is $n$, the index of
$T(\xi)$ is $n$.
\end{enumerate} }
\end{example}

\begin{remark}\hspace{1pt}
{\rm
\begin{enumerate}
\item[(i)] If the index of the Hopf algebra $H$ is greater than 1,
then $H$ is not semisimple by example \ref{index_ex}(i).
\item[(ii)] If $\dim H$ is odd, it
follows from theorem~\ref{Nich89} that the order of the
distinguished group-like element $g$ of $H$ and the order of the
distinguished group-like element $\a$ of $H^*$ are both odd.
Hence, by the formula (\ref{eS}), the order of $S^4$ is also odd.
Therefore, the index of $H$ is odd.
\end{enumerate}
}
\end{remark}
\section{Eigenspace decompositions for Hopf algebras of odd
index}\label{s3}
 In this section, we will only consider those Hopf algebras $H$
 of odd index $n$. Since $r(g)^n=S^{4n}=id_H$, and $S^2$ and $r(g)$
 are commuting operators on $H$, $r(g)$ and $S^2$ are
 simultaneously diagonalizable. Let $\w \in k$ be a primitive $n$th
 root of unity. Then any eigenvalue of $S^2$ is of the form
 $(-1)^a\w^i$ and the eigenvalues of $r(g)$ are of the form $\w^j$.
 Define
$$
H^\w_{a,i,j} = \{ u \in H\,|\, S^2(u) = (-1)^a \w^i u \,, ug =
\w^j u \}\mbox{ for any } (a,i,j) \in \Z_2\times \Z_n \times
\Z_n\,.
$$
 We will simply write $\K_n$ for the group
$\Z_2\times \Z_n \times \Z_n$ and write $H^\w_{(a,i,j)}$ for
$H^\w_{a,i,j}$ for convenience. We then have the decomposition
\begin{equation}\label{eigen_decompose}
H = \bigoplus_{\ba\in\K_n}H^\w_\ba\,.
\end{equation}
Note that $H^\w_\ba$ is not necessarily non-trivial.\\

Since the distinguished group-like element $\a$ of $H^*$ is an
algebra map and $g^n=1$, we have $\a(g)^n=1$. Hence,  $\a(g)$ is a
$n$th root of unity, and so $\a(g) = \w^x$ for some integer $x$.
Using the eigenspace decomposition of $H$ in
(\ref{eigen_decompose}), the diagonalization of the left integral
of $H$ admits an interesting form.

\begin{lem}\label{l4}
Let $H$ be a Hopf algebra over the field $k$ of odd index $n$. Let
$g$ and $\a$ be the distinguished group-like elements of $H$ and
$H^*$, respectively. Suppose that $\L$ is a left integral for $H$.
Then
$$
\Delta(\L) \in \sum_{\ba \in \K_n} H^\w_\ba \otimes
H^\w_{-\ba+\bx}
$$
where $\bx=(0,-x,x)$ and $\a(g)=\w^x$\,.
\end{lem}
\pf Note that
$$
H\otimes H = \bigoplus_{\ba,\bb \in \K_n} H^\w_\ba \otimes
H^\w_\bb\,.
$$
In particular, we can write
$$
\Delta(\L) = \sum_{\ba,\bb \in \K_n} \left(\sum u_\ba \otimes
v_\bb\right)
$$
where $ \sum u_\ba \otimes v_\bb \in H^\w_\ba \otimes H^\w_\bb$.
By~\cite[Proposition~3(d)]{Radf94},
$$
S^2(\L) = \a(g^{-1})\L = \w^{-x}\L\,.
$$
Since $S^2$ is a coalgebra automorphism on $H$, we have
\begin{equation}\label{integral.e1}
\begin{split}
\Delta(\L) &= \sum_{(a,i,j),(b,s,t) \in \K_n}
\left(\sum u_{a,i,j} \otimes v_{b,s,t}\right)\\
 &=
 \sum_{(a,i,j),(b,s,t) \in \K_n}
 \w^{x} S^2\otimes S^2\left(\sum u_{a,i,j} \otimes v_{b,s,t}\right)\\
&=  \sum_{(a,i,j),(b,s,t)\in \K_n} (-1)^{a+b}\w^{x+i+s} \left(\sum
u_{a,i,j} \otimes v_{b,s,t}\right)\,.
\end{split}
\end{equation}
Since $g$ is group-like and $\L g =\a(g)\L=\w^x\L$, we have
\begin{equation}\label{integral.e2}
\begin{split}
\Delta(\L) &= \sum_{(a,i,j),(b,s,t) \in \K_n}
\left(\sum u_{a,i,j} \otimes v_{b,s,t}\right)\\
 &=
 \sum_{(a,i,j),(b,s,t) \in \K_n} \w^{-x}
  r(g)\otimes r(g)\left(\sum u_{a,i,j} \otimes v_{b,s,t}\right)\\
&= \sum_{(a,i,j),(b,s,t) \in \K_n} \w^{-x+j+t}\left(\sum u_{a,i,j}
\otimes v_{b,s,t}\right)\,.
\end{split}
\end{equation}
Thus, if $\sum u_{a,i,j} \otimes v_{b,s,t}\ne 0$, by equations
(\ref{integral.e1}) and (\ref{integral.e2}),
$$
1 = (-1)^{a+b}\w^{x+i+s} \quad \mbox{and}\quad 1=\w^{-x+j+t}\,,
$$
or equivalently,
$$
(b,s,t) = (a,-i,-j)+(0,-x,x)=-(a,i,j)+\bx \,.
$$
Thus,
\begin{equation}\label{normalform}
\Delta(\L) = \sum_{\ba \in \K_n} \left(\sum u_\ba \otimes
v_{-\ba+\bx}\right)\,.
\end{equation}
\hfill\qed
\\

In the sequel, we will call the expression in equation
(\ref{normalform}) the {\em normal form} of $\Delta(\L)$
associated with $\w$. We will simply write $u_\ba \otimes
v_{-\ba+\bx}$ for the sum $\sum u_\ba \otimes v_{-\ba+\bx}$ in the
normal form of $\Delta(\L)$.
\\

The eigenspace decomposition $H = \bigoplus_{\ba \in \K_n}
H^\w_\ba$ is associated with a unique family of projections
$E^\w_\ba$ ($\ba \in \K_n$) from $H$ onto $H^\w_\ba$ such that
\begin{enumerate}
\item $E^\w_\ba \circ E^\w_\bb =0$ for $\ba \ne \bb$ and
\item $\sum_{\ba \in \K_n} E^\w_\ba = id_H$\,.
\end{enumerate}
In particular, $\dim H^\w_{\ba} = \tr(E^\w_\ba)$ for all $\ba \in
\K_n$. By Lemma \ref{l4},
$$
\Delta(\L) = \sum_{\ba\in \K_n} \left(E^\w_\ba \otimes
E^\w_{-\ba+\bx}\right)\Delta(\L)\,
$$
and hence $\left(E^\w_\ba \otimes
E^\w_{-\ba+\bx}\right)\Delta(\L)$ is identical to $\sum u_\ba
\otimes v_{-\ba+\bx}$ in the normal form (\ref{normalform}) of
$\Delta(\L)$. Using the trace formula~\cite[Theorem 1]{Radf90}, we
obtained the following Lemma:

\begin{lem}\label{l5}
Let $H$ be a Hopf algebra over the field $k$ of odd index $n$ and
let $\w \in k$ be a primitive $n$th root of unity. Suppose that
$\L$ is a left integral for $H$ and that $\l$ be a right integral
for $H^*$ such that $\l(\L)=1$. Then
\begin{equation}\label{e6}
\dim   H^\w_\ba =  \l(S(v_{-\ba+\bx} )u_\ba)
\end{equation}
for all $\ba \in \K_n$, where $\sum_{\ba \in \K_n} u_\ba \otimes
v_{-\ba+\bx}$ is the normal form of $\Delta(\L)$ associated with
$\w$.
\end{lem}
\pf  Using the normal form of $\Delta(\L)$ associated with $\w$
and
 \cite[Theorem 1]{Radf90}, for any $\bb \in \K_n$,
\begin{eqnarray*}
\dim H^\w_{\bb} = \tr (E^\w_{\bb}) & = & \sum_{\ba \in \K_n}
\l\left( S(v_{-\ba+\bx})E^\w_{\bb}(u_\ba)\right)\\
& = & \sum_{\ba \in \K_n} \delta_{\ba,\bb}
\l\left( S(v_{-\ba+\bx})u_\ba\right)\\
& = &  \l\left(S(v_{-\bb+\bx})u_{\bb}\right)\,. \qed
\end{eqnarray*}\\

The family of elements $ S(v_{-\ba+\bx} )u_\ba$ appearing in
(\ref{e6}) are in $H^\w_{0,x,0}$. Moreover, if $H$ is
non-semisimple, they satisfy a system of equations.

\begin{lem}\label{l4.5}
Let $H$ be a non-semisimple Hopf algebra over the field $k$ of odd
index $n$ and let $\w \in k$ be a primitive $n$th root of unity.
Then
$$
 \sum_{\substack{(a,i) \in \Z_2 \times \Z_n}}
  (-1)^a\w^{-i}\dim H^\w_{a,i,j} = 0 \mbox{ for } j \in \Z_n\,.
$$
\end{lem}
\pf Let $\L$ be a left integral for $H$ and let $\l$ be a right
integral for $H^*$ such that $\l(\L)=1$. If $H$ is not semisimple,
by \cite[Theorem 4]{Radf94},
$$
\sum S^3(\L_2)\L_1 =0\,.
$$
Hence for any integer $e$,
$$
\sum S^3(\L_2)\L_1g^e =0
$$
Let
$$
h'_\ba =  S^3(v_{-\ba+\bx})u_\ba\quad\mbox{for all } \ba \in
\K_n\,
$$
where $\sum_{\ba \in \K_n} u_\ba \otimes v_{-\ba +\bx}$ is the
normal form of $\Delta(\L)$ associated with $\w$. Then
\begin{eqnarray*}
  0  &= & \sum S^3(\L_2)\L_1g^e \\
  &= &\sum_{(a,i,j)\in  \K_n} h'_{a,i,j} g^e\\
  & = &\sum_{j \in \Z_n} \w^{ej}
    \sum_{(a,i) \in \Z_2 \times \Z_n} h'_{a,i,j}
\end{eqnarray*}
for $e = 0,\dots, n-1$. Since $1, \w,\dots ,\w^{n-1}$ are distinct
elements in $k$, the Vandermonde matrix
$$
\left[\begin{array}{llll}
1&1& \cdots & 1 \\
1&\w & \cdots & \w^{n-1} \\
\vdots & \vdots & & \vdots \\
1 &\w^{n-1} &\cdots & \w^{(n-1)^2}
\end{array}
\right]
$$
is invertible. Therefore,
\begin{equation}\label{temp_eq}
\sum_{(a,i) \in \Z_2 \times \Z_n} h'_{a,i,j}=0
\end{equation}
for $j \in \Z_n$. Notice that
$$
S^3(v_{a,-i-x,-j+x })= (-1)^a\w^{-i-x}S( v_{a,-i-x,-j+x})\,.
$$
Therefore, $h'_{a,i,j} = (-1)^a\w^{-i-x}
S(v_{-a,-i-x,j+x})u_{a,i,j}$ for any $(a,i,j) \in \K_n$. Then
equation (\ref{temp_eq}) becomes
$$
\sum_{(a,i) \in \Z_2 \times \Z_n} (-1)^a\w^{-i}
S(v_{-a,-i-x,j+x})u_{a,i,j}=0
$$
for $j \in \Z_n$. Applying $\l$ to the equation, we  have
$$
\sum_{(a,i) \in \Z_2 \times \Z_n} (-1)^a\w^{-i}
\l\left(S(v_{-a,-i-x,j+x})u_{a,i,j}\right)=0
$$
for all $j \in \Z_n$. Then, the result follows from Lemma
\ref{l5}. \qed
\begin{lem}\label{2.p1}
Let $H$ be a non-semisimple unimodular Hopf algebra over the field
$k$ of odd index $n$. Let $\w \in k$ be a primitive $n$th root of
unity. Then
$$
\sum_{i \in \Z_n}(-1)^a\w^{-i}\dim H^\w_{a,i,l-2i}=0
$$
for $l \in \Z_n$.

\end{lem}
\pf Let $\a$ and $g$ be the distinguished group-like elements of
$H^*$ and $H$, respectively. Since $H$ is unimodular, $\a=\e$  and
hence $\a(g)=1=\w^0$. Let $\L$ be a left integral for $H$ and let
$\l$ be a right integral for $H^*$ such that $\l(\L)=1$. It
follows from Lemma \ref{l4} that the normal form of $\Delta(\L)$
associated with $\w$ is
\begin{equation}\label{sp_normal_form}
\sum_{\ba \in \K_n} u_\ba \otimes v_{-\ba}\,.
\end{equation}
Since $H$ is not semisimple,
$$
0=\e(\L)1=\sum \L_1 S(\L_2)\,.
$$
Thus, we have
\begin{equation}\label{e1.2.p1}
0=\sum_{\ba \in \K_n}  u_\ba S(v_{-\ba})\,.
\end{equation}
Note that, by equation (\ref{eS}) and the unimodularity of $H$,
$$
g^e a = S^{4e}(a)g^e
$$
for any integer $e$ and $a \in H$. Let $\hbar_\ba = u_\ba
S(v_{-\ba})$ for $\ba \in \K_n$. Then,
\begin{equation}
\begin{split}
g^e \hbar_{a,i,j} &=  g^e u_{a,i,j} S(v_{a,-i,-j})\\
&= \w^{e(2i+j)}   u_{a,i,j} S(v_{a,-i,-j}) \\
&=  \w^{e(2i+j)}\hbar_{a,i,j}\,.
\end{split}
\end{equation}
By multiplying $g^e$ on the left in equation (\ref{e1.2.p1}), we
have
\begin{equation}
0=\sum_{(a,i,j)\in \K_n} \w^{e(2i+j)} \hbar_{a,i,j}
 = \sum_{l \in \Z_n} \w^{el} \sum_{(a,i)\in \Z_2 \times \Z_n}
  \hbar_{a,i,l-2i}\,.
\end{equation}
By the same argument used in the proof of Lemma \ref{l4.5},
\begin{equation}\label{e2.2.p1}
\sum_{(a,i) \in \Z_2 \times \Z_n} \hbar_{a,i,l-2i} = 0
\end{equation}
for $l \in \Z_n$. Notice that, by \cite[Theorem 3(a)]{Radf94},
\begin{equation}
\begin{split}
\l(\hbar_{a,i,j}) &=  \l(u_{a,i,j} S(v_{a,-i,-j}))\\
& =  \l(S^3(v_{a,-i,-j})u_{a,i,j})\\
&= (-1)^a\w^{-i}\l(S(v_{a,-i,-j})u_{a,i,j}) \,.
\end{split}
\end{equation}
By Lemma \ref{l5} and equation (\ref{sp_normal_form}),
$$
\l(\hbar_{a,i,j}) = (-1)^a\w^{-i}\dim H^\w_{a,i,j}\,.
$$
Hence, we have
$$
0=\sum_{(a,i)\in \Z_2\times \Z_n}\l(\hbar_{a,i,l-2i})=\sum_{i \in
\Z_n}(-1)^a \w^{-i}\dim H^\w_{a,i,l-2i}\,.
$$
for $l \in \Z_n$. \qed

\section{Arithmetic properties of  Hopf algebras with odd prime index}
\label{s4} In this section, we will study the arithmetic
properties for the Hopf algebras of odd prime index $p$. Let $\w
\in k$ be a primitive $p$th root of unity. The Taft algebra
$T(\w)$ \cite{Taft71} is then a Hopf algebra of this type by
example \ref{index_ex} (ii). The quantum double of $T(\w)$ is a
unimodular Hopf algebra of index $p$ (cf. \cite{RadKa93}).

\begin{lem}\label{l5.5}
Let $H$ be a Hopf algebra of index $p$. Then, for each $j \in
\Z_p$, there exists an integer $d_j$ such that
$$
\dim H^\w_{0,i,j} - \dim H^\w_{1,i,j} = d_j \,.
$$
for any $i \in \Z_p$.
\end{lem}
\pf By Lemma \ref{l4.5}, we have
$$
\sum_{i \in \Z_p} \w^{-i} (\dim H^\w_{0,i,j}- \dim H^\w_{1,i,j})=0
$$
for any $j \in \Z_p$. In particular, $\w^{-1}$ is a root of the
integral polynomial
$$
f_j(x) = \sum_{i=0}^{p-1} (\dim H^\w_{0,i,j} -\dim
H^\w_{1,i,j})\,x^i \,.
$$
Hence, $f_j(x) = d_j \Phi_p(x)$ for some $d_j \in \Q$, where
$\Phi_p(x)=1+x+\cdots+x^{p-1}$ is the irreducible polynomial of
$\w^{-1}$ over $\Q$. Therefore,
$$
\dim H^\w_{0,i,j} -\dim H^\w_{1,i,j} = d_j \,.
$$
Since $\dim H^\w_{0,i,j} -\dim H^\w_{1,i,j}$ is an integer, and so
is $d_j$. \qed

\begin{lem}\label{l6}
Let $H$ be a Hopf algebra of index $p$, where $p$ is an odd prime.
If $H^*$ is not unimodular, then $p \mid \dim H$ and
$$
\sum_{\substack{(a,i)\in \Z_2\times \Z_p}} \dim H^\w_{a,i,j} =
\frac{\dim H}{p}
$$
\end{lem}
\pf Since $H^*$ is not unimodular, $g \ne 1$. Then, $\tr(r(g))=0$
(cf. \cite[Proposition 2.4(d)]{LaRad95}. Moreover, $r(g)^p=id_H$.
Hence, by Lemma \ref{1}, $p|\dim H$ and the eigenspace of $r(g)$
 associated
with the eigenvalue $\w^j$ is of dimension $\frac{\dim H}{p}$ for
any $j \in \Z_p$. Note that
$$
\bigoplus_{(a,i)\in \Z_2\times \Z_p} H^\w_{a,i,j}
$$
is the eigenspace of $r(g)$ associated with $\w^j$. Therefore,
$$
\frac{\dim H}{p} = \dim \left(\bigoplus_{(a,i)\in \Z_2\times \Z_p}
H^\w_{a,i,j} \right) = \sum_{(a,i)\in \Z_2\times \Z_p} \dim
H^\w_{a,i,j}\,. \qed
$$
\\
\begin{lem}\label{l6.5}
Let $H$ be a Hopf algebra of index $p$. If $H^*$ is not unimodular
and $H$ is unimodular, then :
\begin{enumerate}
\item[\rm (i)] There is an integer $d$ such that
$$
\dim H^\w_{0,i,j} - \dim H^\w_{1,i,j} = d \quad \mbox{for
any}\quad i,j \in \Z_p\,.
$$
\item[\rm (ii)] $\tr(S^{2p}) = p^2d$\,.
\end{enumerate}
\end{lem}
\pf (i) By Lemma \ref{2.p1}, for any $l\in \Z_p$,
$$
\sum_{i \in \Z_p} \left(\dim H^\w_{0,i,l-2i} - \dim
H^\w_{1,i,l-2i}\right)\w^{-i} =0
$$
Since $\w^{-1}$ is also a primitive $p$th root of unity in $k$,
there exists an integer $c_l$ such that
\begin{equation}
\dim H^\w_{0,i,l-2i} - \dim H^\w_{1,i,l-2i} = c_l
\end{equation}
for $i \in \Z_p$. By Lemma \ref{l5.5}, for any $i,l \in \Z_p$,
\begin{equation}
c_l=\dim H^\w_{0,i,l-2i} - \dim H^\w_{1,i,l-2i} = d_{l-2i}\,.
\end{equation}
Since $2$ and $p$ are relative prime, $l,l-2,\dots,l-2(p-1)$ is a
complete set of representatives of $\Z_p$. Therefore,
$$
d_j = c_l = d \quad \mbox{ for any }j,l \in \Z_p\,.
$$
(ii) Since $p$ is odd,
\begin{eqnarray*}
\tr(S^{2p})&=&
\sum_{i,j \in \Z_p} \dim H^\w_{0,i,j} - \dim H^\w_{1,i,j}\\
&=& \sum_{i,j \in \Z_p} d \\
&=&p^2d\,.\qed
\end{eqnarray*}
\section{Hopf algebras of dimension $pq$}\label{s5}
In this section, we will consider the Hopf algebras $H$ of
dimension $p q$ where both $p \le q$ are odd primes. In
particular, we prove that if $H$ is not semisimple and $\dim H
=p^2$, then $H$ is isomorphic to a Taft algebra. By \cite[Theorem
2]{Mas96}, any Hopf algebra over $k$ of dimension $p^2$ is either
a group algebra or a Taft algebra. We begin the section with the
following lemma.

\begin{lem}\label{l5.1}
Let $p,q$ be two distinct prime numbers. Then there is no Hopf
algebra $H$ of dimension $p q$ such that $|G(H)|=p$ and
$|G(H^*)|=q$.
\end{lem}
\pf Suppose there is a Hopf algebra $H$ of dimension $pq$ such
that $|G(H)|=p$ and $|G(H^*)|=q$. Let $g \in G(H)$ and $\a \in
G(H^*)$ such that $o(g)=p$ and $o(\a)=q$. Then,
$$
\a(g)^p=\a(g^p) = \a(1)=1
$$
and
$$
1 =\e(g) = \a^q(g) =\a(g)^q\,.
$$
Therefore, $o(\a(g))=1$ and so $\a(g)=1$. Since $k[G(H^*)]$ is a
Hopf subalgebra of $H^*$, $k[G(H^*)]^\perp$ is a Hopf ideal of
$H$. Let $B^+$ be the augmentation ideal of $k[G(H)]$. Then $B^+H
= (g-1)H$ and
$$
\a^i((g-1)h)=(\a(g)^i-1)\a^i(h)=0\quad \mbox{for}\quad h \in H, i
\in \Z\,.
$$
Therefore,
$$
B^+H \C k[G(H^*)]^\perp\,.
$$
It follows from \cite[Theorem 2.4 (2a)]{Schn92} that $\dim H/B^+H
= q$. Thus,
$$
\dim B^+H = pq - q = k[G(H^*)]^\perp
$$
and hence,
$$
B^+H = k[G(H^*)]^\perp\,.
$$
Therefore, $H/B^+H$ is isomorphic to $k[G(H^*)]^*$ as Hopf
algebras. In particular, $H/B^+H$ is semisimple. Let $\L$ be a
non-zero left integral of $H$ and $\L'$ a non-zero right integral
of $k[G(H)]$. Since $char k =0$, $\e(\L') \ne 0$ and hence,
$\L'\L=\e(\L')\L \ne 0$. Therefore, $\L \not\in B^+H$ and so $\L+
B^+H$ is a non-zero left integral in $H/B^+H$. Since $H/B^+H$ is
semisimple, $\e(\L)=\e(\L+B^+H) \ne 0$. Hence, $H$ is semisimple.
By \cite{EG99}, $H$ is trivial and so $|G(H)|=pq$, a
contradiction. \hfill$\blacksquare$

\begin{prop}\label{index_prop}
Let $H$ be a  non-semisimple Hopf algebra of dimension $p q$ where
$p \le q$ are odd primes. Then
\begin{enumerate}
\item[\rm (i)] the order of $S^4$ is $p$ and
\item[\rm (ii)] $H$ is of index $p$.
\end{enumerate}
\end{prop}
\pf (i) Since $H$ is not semisimple and $\dim H$ is odd, by
\cite[Theorem 2.1]{LaRad95} or \cite[Lema 2.5]{AnSc98}, $S^4 \ne
id_H$ and $H$, $H^*$ cannot both be unimodular. Let $g$ be the
distinguished group-like element of $H$ and let $\a$ the
distinguished group-like element of $H^*$. Then,
 $o(\a)<pq$ and $o(g) <  pq$, for otherwise, $H$ is isomorphic
to a group algebra  which is semisimple. By Lemma \ref{l5.1},
\begin{equation}\label{e5.1}
\lcm(o(g),o(\a)) = p \mbox{ or } q\,.
\end{equation}
By the equation (\ref{eS}) and (\ref{e5.1}), the order of $S^4$ is
either $p$ or $q$. If $p=q$, order of $S^4$ and the index of $H$
are obviously equal to $p$. We now assume $q >p$. We consider the
following cases:\\
{\em Case (a):  $H^*$ is not unimodular.} Suppose that the order
of $S^4$ is  $q$. By equation (\ref{eS}), $q \mid
\lcm(o(g),o(\a))$. Therefore, $\lcm(o(g),o(\a))=q$  and hence
$o(g)=1$ or $q$. Thus, the index of $H$ is also $q$. Let $\w\in k$
be a $q$th primitive root of unity. By Lemma \ref{l5.5}, for each
$j \in \Z_q$ there is an integer $d_j$ such that
\begin{equation}\label{e5.2}
\dim H^\w_{0,i,j} - \dim H^\w_{1,i,j} = d_j \quad \mbox{for all }
i \in \Z_q\,.
\end{equation}
Let $X_{i,j} = \min(\dim H^\w_{0,i,j}, \dim H^\w_{1,i,j})$. Then,
$$
\dim H^\w_{0,i,j}+ \dim H^\w_{1,i,j} = 2X_{i,j} + |d_j|
$$
and so
\begin{equation}\label{e5.3}
\sum_{(a,i) \in \Z_2 \times \Z_q } H^\w_{a,i,j} = \sum_{i \in \Z_q
} 2X_{i,j} + q|d_j|
\end{equation}
for each $j \in \Z_q$. It follows from Lemma \ref{l6} that
\begin{equation}\label{e5.4}
\sum_{i \in \Z_q} 2X_{i,j} + q|d_j| = p\,.
\end{equation}
Since $p$ odd, by (\ref{e5.4}), $|d_j|$ must be odd. However, the
left hand side of (\ref{e5.4}) is then strictly greater than $p$,
a contradiction! Therefore, $o(S^4)=p$.\\
{\em Case (b):  $H^*$ is unimodular}. Then $H^{**} \cong H$ is not
unimodular. By Theorem \ref{th1.1}, $H^*$ is not semisimple and
$\dim H^*  =pq$. It follows from Case (a) that the order of
${S^*}^4$
is $p$. Since $o(S^4) = o({S^*}^4)$. Therefore, $o(S^4)=p$. \\
(ii) Let $n$ be the index of $H$. Then, by (\ref{eqdiv}), $n \mid
\lcm(o(g),o(\a))$ and $o(S^4) \mid n$. Since $o(S^4)=p$ and
$\lcm(o(g),o(\a)) = p \mbox{ or } q$, we have $n=p$.  \qed

\begin{lem}\label{biproduct}
Let $H$ be a Hopf algebra over $k$ such that both the
distinguished group-like elements of $H$ and $H^*$ are of order
$p$ where $p$ is an odd prime. Then, $\tr(S^{2p})=p^2 d$ for some
integer $d$.
\end{lem}
\pf Let $g$ and $\a$ be the distinguished group-like elements of
$H$ and $H^*$ respectively. Let $B$ be the group algebra $k[g]$.
It follows from the arguments in the proof of \cite[Theorem
A]{AnSc98} that there is an Hopf algebra map $\pi: H \map B$ such
that $\pi \gamma = id_B$ where $\gamma : B \map H$ is the
inclusion map. Therefore, $H$ isomorphic to the biproduct $R
\times B$ as Hopf algebras where
$$
R=H^{co B}=\{h \in H| (id
\otimes \pi)\Delta(h)=h \otimes 1\}
$$ (cf. \cite{Radf85}). It is
shown in \cite[section 4 ]{AnSc98} that $R$ is invariant under
$S^2$. Moreover, in the identification $H \cong R \otimes H$ given
by multiplication, one has
\begin{equation}\label{S_square}
S^2 = T \otimes id_B\,.
\end{equation}
Since $H$ is not unimodular, $H$ is not semisimple and hence
$\tr(S^2)=0$. By equation (\ref{S_square}), $\tr(S^2)=\tr(T)p$.
Therefore, $\tr(T)=0$. Moreover, $T^{2p}=id_R$ as $S^{4p}=id_H$ by
equation (\ref{eS}). Hence, by Lemma \ref{1}, $\tr(T^p)=pd$ for
some integer $d$. Since $ S^{2p}=T^p \otimes id_B$, we have
$$
\tr(S^{2p}) = \tr(T^p) \tr(id_B) = p^2d. \qed
$$

\begin{thm}\label{main_pq}
Let $H$ be a non-semisimple Hopf algebra of dimension $p q$ where
$p \le q$ are odd primes. Then $\tr(S^{2p}) = p^2 d$ for some  odd
integer $d$\,.
\end{thm}
\pf By Proposition \ref{index_prop}, $S^{4p}=id_H$. Let
$$
H_\pm = \{h \in H| S^{2p}(h) = \pm h\}\,.
$$
Then,
\begin{eqnarray*}
\dim H_+ - \dim H_- & = & \tr(S^{2p})\\
\mbox{and}\quad    \dim H_+ + \dim H_- &=& pq\,.
\end{eqnarray*}
Since $pq$ is odd, $\tr(S^{2p})$  is also an odd integer. Thus, if
$\tr(S^{2p})=p^2 d$, then $d$ must be an odd integer. Therefore,
it suffices to show that $\tr(S^{2p})=p^2 d$ for some integer $d$.
Since $H$ is not semisimple, by Theorem \ref{th1.1}, $H^*$ is also
not semisimple. By Proposition \ref{index_prop}, the indexes of
$H$ and $H^*$ are both $p$. Since $\dim H$ is odd, by
\cite[Theorem 2.2]{LaRad95}, not both of $H$ and $H^*$ are
unimodular. We then have the following three cases:\\
(i) If $H$ is unimodular and $H^*$ is not unimodular,  the result
follows from Lemma \ref{l6.5}.
\\
(ii) If $H$ is not unimodular and $H^*$ is unimodular, by Lemma
\ref{l6.5}, $\tr({S^*}^{2p})=p^2d$ for some odd integer $d$. The
result follows from $\tr({S^*}^{2p})= \tr(S^{2p})$.
\\
(iii) If both $H$ and $H^*$ are not unimodular, by Lemma
\ref{l5.1} and Proposition \ref{index_prop}, the orders of the
distinguished group-like elements of $H$ and $H^*$ are both equal
to $p$. Thus, by
Lemma \ref{biproduct}, $\tr(S^{2p})=p^2 d$. \qed\\

As a consequence of the above theorem, we prove that any Hopf
algebra of dimension $p^2$ is either a group algebra or a Taft
algebra (see example \ref{index_ex}(ii)).

\begin{thm}\label{main_pp}
Let $H$ be a Hopf algebra over $k$ of dimension $p^2$ where $p$ is
any prime number. Then, $H$ is isomorphic to one of the following
Hopf algebras:
\begin{enumerate}
\item[\rm (a)] $k[\Z_{p^2}]$ ;
\item[\rm (b)] $k[\Z_p \times \Z_p]$;
\item[\rm (c)] $T(\w)$, $\w\in k$ a primitive $p$th of unity.
\end{enumerate}
\end{thm}
\pf If $H$ is semisimple, it follows from \cite[Theorem 2]{Mas96}
that $H$ isomorphic to  $k[\Z_{p^2}]$ or $k[\Z_p \times \Z_p]$. It
is also shown in \cite{Kapl75} that if $H$ is a non-semisimple
Hopf algebra of dimension 4, then
 $H$ isomorphic to the Taft algebra $T(1)$ or $T(-1)$.
We may now assume $H$ is not
 semisimple and $p$ is odd. Let $S$ be the antipode of $H$.
  By Proposition \ref{index_prop},
 $S^{4p}=id_H$ and so $S^{2p}$ is
 diagonalizable and  the possible eigenvalues of $S^{2p}$
 are $\pm 1$. Suppose $S^{2p} \ne id_H$. Then, $\tr(S^{2p})$ is
 an integer such that
$$
-p^2  \le \tr(S^{2p}) < p^2\,.
$$
By Theorem
 \ref{main_pq},
 $$
 \tr(S^{2p}) =p^2 d
 $$
 for some odd integer $d$. Therefore, $\tr(S^{2p})=- p^2$ and hence
 $S^{2p}=-id_H$. However, this is not possible since $S^{2p}(1_H)= 1_H$.
 Therefore, $S^{2p} = id_H$.
  By Proposition \ref{index_prop},
  the order of $S^4$ is $p$, and so is the order $S^2$. It follows
  from \cite[Theorem A(ii)]{AnSc98} that $H$ is isomorphic to a
  Taft algebra of dimension $p^2$. Hence, $H \cong T(\w)$
  for some primitive $p$th root of unity, $\w \in k$.
 \qed
\begin{center}
{\bf Acknowledgement}
\end{center}

The author would like to thank Susan Montgomery for bringing his
attention to the question on the Hopf algebras of dimension $p^2$.

\bibliographystyle{amsalpha}
\providecommand{\bysame}{\leavevmode\hbox
to3em{\hrulefill}\thinspace}
\providecommand{\MR}{\relax\ifhmode\unskip\space\fi MR }
\providecommand{\MRhref}[2]{%
  \href{http://www.ams.org/mathscinet-getitem?mr=#1}{#2}
} \providecommand{\href}[2]{#2}

\vspace{0.5cm}
\begin{quote}
{\bf e-mail:} {\tt rng@towson.edu}
\end{quote}
\end{document}